\documentclass[12pt]{amsart}
\usepackage{amssymb}
\usepackage[all]{xy}

\newcommand{\issuenumber}{10}
\newcommand{\issuemonth}{September}
\newcommand{\issueyear}{2004}

\newtheorem{thm}{Theorem}[section]
\newtheorem{prob}[thm]{Problem}

\newtheorem{lem}[thm]{Lemma}

\newtheorem{issue}{Issue}

\theoremstyle{definition}
\newtheorem{defn}[thm]{Definition}

\theoremstyle{remark}

      \newenvironment{changemargin}[2]{\begin{list}{}{
         \setlength{\topsep}{0pt}\setlength{\leftmargin}{0pt}
         \setlength{\rightmargin}{0pt}
         \setlength{\listparindent}{\parindent}
         \setlength{\itemindent}{\parindent}
         \setlength{\parsep}{0pt plus 1pt}
         \addtolength{\leftmargin}{#1}\addtolength{\rightmargin}{#2}
         }\item }{\end{list}}

\newcommand{\fb}[1]{\fbox{$#1$}}
\newcommand{\nop}{$\times$}
\newcommand{\fbn}{\!\!\fbox{\!\nop\!}\!\!}
\newcommand{\yup}{\checkmark}

\newcommand{\mbq}{\mb{?}}
\newcommand{\mb}[1]{{\mbox{\textbf{#1}}}}
\newcommand{\smb}[1]{{\!\!\mb{#1}\!\!}}
\newcommand{\x}{\times}

\newcommand{\Cantor}{{{}^\N\{0,1\}}}
\newcommand{\oo}{\infty}

\newcommand{\sr}[2]{{\txt{$#1$\\$#2$}}}

\renewcommand{\b}{\mathfrak{b}}
\newcommand{\g}{\mathfrak{g}}

\renewcommand{\d}{\mathfrak{d}}
\newcommand{\p}{\mathfrak{p}}
\newcommand{\s}{\mathfrak{s}}

\newcommand{\NON}{{\mathsf   {NON}}}
\newcommand{\COF}{{\mathsf   {COF}}}

\newcommand{\upannouncement}[1]{[\S\ref{#1} above]}
\newcommand{\dnannouncement}[1]{[\S\ref{#1} below]}

\newcommand{\cA}{\mathcal{A}}
\newcommand{\M}{\mathcal{M}}

\newcommand{\cov}{\mathsf{cov}}

\newcommand{\cf}{\mathsf{cf}}
\newcommand{\non}{\mathsf{non}}

\newcommand{\CH}{the Continuum Hypothesis}

\newcommand{\fo}{\mathfrak{od}}

\newcommand{\impl}{\rightarrow}
\newcommand{\arrays}{{{\{0,1\}}^{\N\x\N}}}
\newcommand{\w}{\omega}
\renewcommand{\b}{\mathfrak{b}}
\renewcommand{\t}{\mathfrak{t}}
\newcommand{\h}{\mathfrak{h}}

\renewcommand{\split}{\mathsf{Split}}
\newcommand{\bq}{\begin{quote}}
\newcommand{\eq}{\end{quote}}
\renewcommand{\O}{\mathcal{O}}
\newcommand{\B}{\mathcal{B}}
\newcommand{\BG}{\B_\Gamma}

\newcommand{\BT}{\B_\Tau}

\newcommand{\BO}{\B_\Omega}
\newcommand{\CG}{C_\Gamma}

\newcommand{\CT}{C_\Tau}

\newcommand{\CO}{C_\Omega}
\newcommand{\sone}{\mathsf{S}_1}    \newcommand{\sfin}{\mathsf{S}_{fin}}
\newcommand{\ufin}{\mathsf{U}_{fin}}

\newcommand{\nin}{\not\in}

\newcommand{\cU}{\mathcal{U}}

\newcommand{\NN}{{{}^{\naturals}\naturals}}
\newcommand{\naturals}{{\mathbb N}}
\newcommand{\N}{\naturals}
\newcommand{\sbst}{\subseteq}
\newcommand{\by}[2]{\par\hfill\emph{#1}, #2}
\newcommand{\nby}[1]{\par\hfill\emph{#1}}
\newcommand{\Tau}{\mathrm{T}}
\newcommand{\CE}{\textsc{CE}}

\newcommand{\be}{\begin{enumerate}}
\newcommand{\ee}{\end{enumerate}}
\newcommand{\bi}{\begin{itemize}}
\newcommand{\ei}{\end{itemize}}
\renewcommand{\i}{\item}

\newcommand{\general}{\small\vfill\par\noindent\hrulefill\par
\noindent\textbf{Previous issues.} The first issues of this
bulletin, which contain general information (first issue), basic
definitions, research announcements, and open problems (all
issues) are available online, on \arx{math.GN/$x$}, where $x$ is
\texttt{0301011}, \texttt{0302062}, \texttt{0303057},
\texttt{0304087}, \texttt{0305367}, \texttt{0312140},
\texttt{0401155}, \texttt{0403369}, and \texttt{0406411},
respectively, for issues number $1$ to $9$.\\[0.1cm]
\textbf{Contributions.}
Please submit your contributions (announcements, discussions, and open problems)
by e-mailing us. It is preferred to write them
in \LaTeX{}.
The authors are urged to use as standard notation as possible, or otherwise give
the definitions or a reference to where the notation is explained.
Contributions to this bulletin would not require any transfer of copyright,
and material presented here can be published elsewhere.\\[0.1cm]
\textbf{Subscription.}
To receive this bulletin (free) to your
e-mailbox, e-mail us:\\
{tsaban@math.huji.ac.il}
}

\newcommand{\nArxPaper}[5]{\subsection{#2}{#4}\par\hfill{\arx{#1}}\par\hfill\emph{#3}}

\newcommand{\nAMSPaper}[5]{\subsection{#2}{#4}\par\hfill{\texttt{#1}}\par\hfill\emph{#3}}
\newcommand{\SPMBul}{\textbf{$\mathcal{SPM}$ Bulletin}}

\newcommand{\arx}[1]{\texttt{http://arxiv.org/abs/#1}}
\newcommand{\url}[1]{\bq\texttt{#1}\eq}
\newcommand{\online}[1]{The paper is available online at \url{#1}}

\newcommand{\probmonth}{\emph{Problem of the month}}

\title[$\mathcal{SPM}$ Bulletin \textbf{\issuenumber} (\issuemonth{} \issueyear)]{%
$\mathcal{SPM}$ Bulletin\\[0.5cm]
Issue number \issuenumber: \issuemonth{} \issueyear{} \CE{}}

\begin{document}
\maketitle

\tableofcontents

\section{Editor's note}

We are glad to announce the solution of $4+21+\frac{1}{2}$ (!)
problems posed in earlier issues of the \SPMBul{};
the ``$\frac{1}{2}$'' standing for a ``consistently yes'' answer
of Zdomsky to the last issue's \probmonth{}.

\medskip

Contributions to the next issue are, as always, welcome.

\medskip

\by{Boaz Tsaban}{tsaban@math.huji.ac.il}

\hfill \texttt{http://www.cs.biu.ac.il/\~{}tsaban}

\section{Research announcements}

\newcommand{\mad}{\operatorname{MAD}}
\nAMSPaper{http://www.ams.org/journal-getitem?pii=S0002-9939-04-07580-X}
{On two problems of Erd\H os and Hechler: New methods in singular Madness}
{Menahem Kojman, Wieslaw Kubis, and Saharon Shelah}
{For an infinite cardinal $\mu$, $\mad(\mu)$
denotes the set of all cardinalities of \emph{nontrivial maximal
almost disjoint families} over $\mu$.

Erd\H os and Hechler proved in \cite{EH} the consistency of
$\mu\in \mad(\mu)$  for a singular cardinal $\mu$ and asked if it
was ever possible for a singular $\mu$ that $\mu\notin
\mad(\mu)$, and also whether $2^{\cf\mu}<\mu \Longrightarrow
\mu\in \mad(\mu)$ for every singular cardinal $\mu$.

We introduce a new method for controlling $\mad (\mu)$ for a
singular $\mu$ and, among other new results about the structure of
$\mad(\mu)$ for singular $\mu$,  settle both problems
affirmatively.}

\nArxPaper{math.LO/0406530}
{Almost isometric embeddings of metric spaces}
{Menachem Kojman and Saharon Shelah}
{We investigate a relations of almost isometric embedding and almost isometry
between metric spaces and prove that with respect to these relations:
\be
\i There is a countable universal metric space.
\i There may exist fewer than continuum separable metric spaces on $\aleph_1$
so that every separable metric space is almost isometrically embedded into one
of them when the continuum hypothesis fails.
\i There is no collection of fewer than continuum metric spaces of
cardinality $\aleph_2$ so that every ultra-metric space of cardinality $\aleph_2$ is
almost isometrically embedded into one of them if $\aleph_2<2^{\aleph_0}$.
\ee
We also prove that various spaces $X$ satisfy that if a space $Y$ is almost
isometric to $X$ than $Y$ is isometric to $X$.}

\nArxPaper{math.LO/0407405}
{Half of an inseparable pair}
{Arnold W.\ Miller}
{A classical theorem of Luzin is that the separation principle holds for the
$\Pi^0_\alpha$ sets but fails for the $\Sigma^0_\alpha$ sets. We show that for every
$\Sigma^0_\alpha$ set $A$ which is not $\Pi^0_\alpha$ there exists a $\Sigma^0_\alpha$ set
$B$ which is disjoint from A but cannot be separated from $A$ by a $\Delta^0_\alpha$ set
$C$. Assuming $\Pi^1_1$-determinacy it follows from a theorem of Steel that a similar
result holds for $\Pi^1_1$ sets. On the other hand assuming $V=L$ there is a proper
$\Pi^1_1$ set which is not half of a Borel inseparable pair. These results answer
questions raised by F.\ Dashiell.}

\nArxPaper{math.LO/0407487}
{Covering the Baire space by families which are not finitely dominating}
{Heike Mildenberger, Saharon Shelah, and Boaz Tsaban}
{It is consistent (relative to ZFC) that each union of $\max\{\b,\g\}$ many families
in the Baire space which are not finitely dominating is not dominating.
In particular:
\be
\i It is consistent that each union of $\max\{\b,\g\}$ many sets with the Scheepers
property $\ufin(\O,\Omega)$ has Menger's property $\ufin(\O,\O)$.
\i It is consistent that for each nonprincipal ultrafilter $\cU$, the
cofinality of the reduced ultrapower $\NN/\cU$ is greater than $\max\{\b,\g\}$.
\ee
The model is constructed by oracle chain condition forcing, to which we give a
self-contained introduction.}

\nArxPaper{math.GN/0408039}
{Compact Scattered Spaces in Forcing Extensions}
{Kenneth Kunen}
{We consider the cardinal sequences of compact scattered spaces in models
where \CH{} is false. We describe a number of models where the continuum is
$\aleph_2$ in which no such space can have $\aleph_2$ countable levels.}

\nArxPaper{math.GN/0408134}
{A note on an example by van Mill}
{K.\ P.\ Hart and G.\ J.\ Ridderbos}
{Improving on an earlier example by J.\ van Mill, we prove that there exists a
zero-dimensional compact space of countable $\pi$-weight and uncountable character
which is homogeneous under $MA+\lnot CH$, but not under $CH$.}

\nAMSPaper{http://www.ams.org/journal-getitem?pii=S0002-9939-04-07376-9}
{Spaces on which every pointwise convergent series of continuous functions converges pseudo-normally}
{Lev Bukovsky and Krzysztof Ciesielski}
{A topological space $X$ is a $\Sigma\Sigma^*$-space
provided that, for every sequence $\langle
 f_n\rangle_{n=0}^\infty$ of continuous functions
from $X$ to $\mathbb{R}$, if the series $\sum_{n=0}^\infty|f_n|$ converges
pointwise, then it converges pseudo-normally.
We show that every regular Lindel\"of
$\Sigma\Sigma^*$-space has the Rothberger property. We also construct,
under the continuum hypothesis, a $\Sigma\Sigma^*$-subset of $\mathbb{R}$ of
cardinality continuum.
}

\subsection{A semifilter approach to selection principles}\label{addsplit}
In this paper we develop the semifilter approach to the classical Menger and Hurewicz
properties and show that the small cardinal $\mathfrak g $ is a
 lower bound of the additivity number of the $\sigma$-ideal generated by
 Menger subspaces of the Baire space, and under $\mathfrak u<\mathfrak g$
every subset $X$ of the real line with the property $\mathrm{Split}(\Lambda,\Lambda)$ is Hurewicz,
and thus it is consistent
with ZFC that the property $\mathrm{Split}(\Lambda,\Lambda)$ is preserved by unions of less than $\mathfrak b$
subsets of the real line.
\nby{Lubomyr Zdomsky}

\label{cot0}
\nArxPaper{math.GN/0409068}
{The combinatorics of $\tau$-covers}
{Heike Mildenberger, Saharon Shelah, and Boaz Tsaban}
{We solve four out of the six open problems concerning critical
cardinalities of topological diagonalization properties involving
$\tau$-covers, show that the remaining two cardinals are
equal, and give a consistency result concerning this remaining cardinal.
Consequently, $21$ open problems concerning potential
implications between these properties are settled. We also give
structural results based on the combinatorial techniques.\\
For further details see \dnannouncement{cot}.
}

\section{The combinatorics of $\tau$-covers: $4+21$ problems solved}\label{cot}
In Section 4 of the fourth issue of the \SPMBul{},
a discussion is made of the Scheepers diagram of selection hypotheses,
when extended by allowing $\tau$-covers.
$6$ critical cardinalities in that diagram were unknown.
Moreover, $76$ potential implications among properties in this
diagram were unsettled.

These remaining cardinalities are addressed in \upannouncement{cot0}.
We give here one example of that treatment, and quote the main results.
Everything in the remainder of this subsection is quoted, without further notice,
from the paper \upannouncement{cot0}.

Let $\CG$, $\CT$, and $\CO$ denote the collections of
\emph{clopen} $\gamma$-covers, $\tau$-covers, and $\w$-covers of $X$, respectively.
Recall that, since we are dealing with sets of reals, we may assume that all open covers
are countable. Restricting attention to countable covers,
we have the following, where an arrow denotes inclusion:
$$\begin{matrix}
\BG      & \impl & \BT      & \impl & \BO      & \impl & \B      \\
\uparrow &       & \uparrow &       & \uparrow &       & \uparrow \\
\Gamma   & \impl & \Tau     & \impl & \Omega   & \impl & \O  \\
\uparrow &       & \uparrow &       & \uparrow &       & \uparrow \\
\CG      & \impl & \CT      & \impl & \CO      & \impl & C
\end{matrix}$$
As each of the properties $\Pi(\cdot\ ,\cdot)$, $\Pi\in\{\sone,\sfin,\ufin\}$,
is monotonic in its first variable, we have that for each $x,y\in\{\Gamma, \Tau, \Omega,\O\}$,
$$\Pi(\B_x,\B_y)\impl\Pi(x,y)\impl\Pi(C_x,C_y)$$
(here $C_\O:=C$ and $\B_\O:=\B$).
Consequently,
$$\non(\Pi(\B_x,\B_y))\le\non(\Pi(x,y))\le\non(\Pi(C_x,C_y)).$$

\begin{defn}\label{diagble}
We use the short notation $\forall^\oo$ for ``for all but finitely many''
and $\exists^\oo$ for ``there exist infinitely many''.
\be
\i $A\in \arrays$ is a \emph{$\gamma$-array}
if $(\forall n)(\forall^\oo m)\ A(n,m)=1$.
\i $\cA\sbst\arrays$ is a $\gamma$-family if each $A\in\cA$ is a
$\gamma$-array.
\i A family $\cA\sbst\arrays$ is \emph{finitely $\tau$-diagonalizable} if
there exist finite (possibly empty) subsets $F_n\sbst\N$, $n\in\N$, such
that:
\be
\i For each $A\in\cA$: $(\exists^\oo n)(\exists m\in F_n)\ A(n,m)=1$;
\i For each $A,B\in\cA$:\\
\begin{tabular}{ll}
Either  & $(\forall^\oo n)(\forall m\in F_n)\ A(n,m)\le B(n,m)$,\\
or      & $(\forall^\oo n)(\forall m\in F_n)\ B(n,m)\le A(n,m)$.
\end{tabular}
\ee
\ee
\end{defn}

Translating the notions of covers into corresponding combinatorial notions,
one obtains the following.
(Notice that $\arrays$ is topologically the same as the Cantor space $\Cantor$.)

\begin{thm}\label{charSfinGT}
For a set of reals $X$, the following are equivalent:
\be
\i $X$ satisfies $\sfin(\BG,\BT)$; and
\i For each Borel function $\Psi:X\to\arrays$, if $\Psi[X]$ is a $\gamma$-family,
then it is finitely $\tau$-diagonalizable.
\ee
The corresponding assertion for $\sfin(\CG,\CT)$
holds when ``Borel'' is replaced by ``continuous''.
\end{thm}

Then the following is proved.
\begin{lem}\label{dgbl=b}
The minimal cardinality of a $\gamma$-family which is not finitely $\tau$-diagonaliz\-able is $\b$.
\end{lem}
Having Theorem \ref{charSfinGT} and Lemma \ref{dgbl=b}, we get that
$$\b=\non(\sfin(\BG,\BT))\le\non(\sfin(\Gamma,\Tau))\le\non(\sfin(\CG,\CT))=\b,$$
and therefore $\non(\sfin(\Gamma,\Tau))=\b$.
It follows that $\non(\sone(\Gamma,\Tau))=\b$, and using a similar approach,
it is proved that $\non(\sone(\Tau,\Tau))=\t$, and $\non(\sfin(\Tau,\Tau))=\min\{\s,\b\}$.

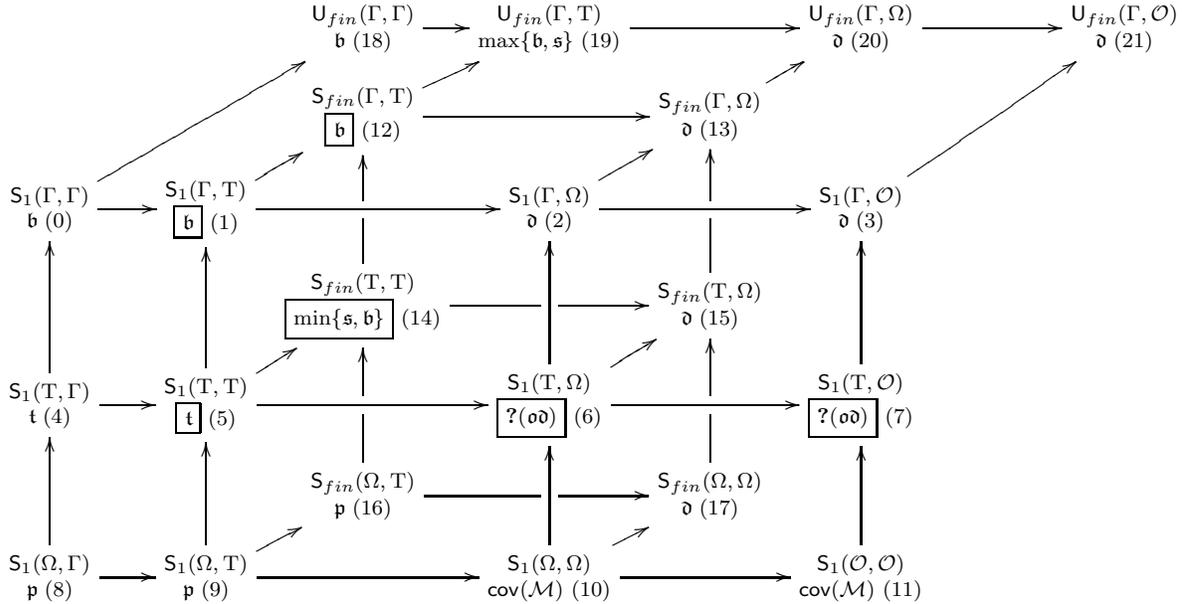
\begin{figure}[!ht]
\renewcommand{\sr}[2]{{\txt{$#1$\\$#2$}}}
{\tiny
\begin{changemargin}{-3cm}{-3cm}
\begin{center}
$\xymatrix@C=7pt@R=6pt{
%1
&
&
& \sr{\ufin(\Gamma,\Gamma)}{\b~~ (18)}\ar[r]
& \sr{\ufin(\Gamma,\Tau)}{\max\{\b,\s\}~~ (19)}\ar[rr]
&
& \sr{\ufin(\Gamma,\Omega)}{\d~~ (20)}\ar[rrrr]
&
&
&
& \sr{\ufin(\Gamma,\O)}{\d~~ (21)}
\\
%2
&
&
& \sr{\sfin(\Gamma,\Tau)}{\fb{\b}~~ (12)}\ar[rr]\ar[ur]
&
& \sr{\sfin(\Gamma,\Omega)}{\d~~ (13)}\ar[ur]
\\
%3
\sr{\sone(\Gamma,\Gamma)}{\b~~ (0)}\ar[uurrr]\ar[rr]
&
& \sr{\sone(\Gamma,\Tau)}{\fb{\b}~~ (1)}\ar[ur]\ar[rr]
&
& \sr{\sone(\Gamma,\Omega)}{\d~~ (2)}\ar[ur]\ar[rr]
&
& \sr{\sone(\Gamma,\O)}{\d~~ (3)}\ar[uurrrr]
\\
%4
&
&
& \sr{\sfin(\Tau,\Tau)}{\fb{\min\{\s,\b\}}~~ (14)}\ar'[r][rr]\ar'[u][uu]
&
& \sr{\sfin(\Tau,\Omega)}{\d~~ (15)}\ar'[u][uu]
\\
\sr{\sone(\Tau,\Gamma)}{\t~~ (4)}\ar[rr]\ar[uu]
&
& \sr{\sone(\Tau,\Tau)}{\fb{\t}~~ (5)}\ar[uu]\ar[ur]\ar[rr]
&
& \sr{\sone(\Tau,\Omega)}{\fbox{\textbf{?}$(\fo)$}~~ (6)}\ar[uu]\ar[ur]\ar[rr]
&
& \sr{\sone(\Tau,\O)}{\fbox{\textbf{?}$(\fo)$}~~ (7)}\ar[uu]
\\
&
&
& \sr{\sfin(\Omega,\Tau)}{\p~~ (16)}\ar'[u][uu]\ar'[r][rr]
&
& \sr{\sfin(\Omega,\Omega)}{\d~~ (17)}\ar'[u][uu]
\\
\sr{\sone(\Omega,\Gamma)}{\p~~ (8)}\ar[uu]\ar[rr]
&
& \sr{\sone(\Omega,\Tau)}{\p~~ (9)}\ar[uu]\ar[ur]\ar[rr]
&
& \sr{\sone(\Omega,\Omega)}{\cov(\M)~~ (10)}\ar[uu]\ar[ur]\ar[rr]
&
& \sr{\sone(\O,\O)}{\cov(\M)~~ (11)}\ar[uu]
}$
\end{center}
\end{changemargin}
}

\caption{The Scheepers diagram, enhanced with $\tau$-covers}\label{tauSch}
\end{figure}

It is not difficult to see that $\non(\sfin(\Tau,\Omega))=\non(\sfin(\Tau,\O))$,
call this joint cardinal $\fo$, the \emph{$o$-diagonalization number}\index{$o$-diagonalization number};
the reason for this to be explained soon.
The surviving properties (in the open case)
appear in Figure \ref{tauSch}, with their critical cardinalities,
and serial numbers (for later reference).
The newly found cardinalities are framed.

By Figure \ref{tauSch},
$$\cov(\M)=\non(\sone(\O,\O))\le\non(\sone(\Tau,\O))\le\non(\sone(\Gamma,\O))=\d,$$
thus $\cov(\M)\le\fo\le\d$.

\begin{defn}\label{odiagbl}
A $\tau$-family $\cA$ is \emph{$o$-diagonalizable} if
there exists a function $g:\N\to\N$, such that:
$$(\forall A\in\cA)(\exists n)\ A(n,g(n))=1.$$
\end{defn}

As in Theorem \ref{charSfinGT}, $\sone(\BT,\B)$ and $\sone(\CT,C)$
have a natural combinatorial characterization.
This characterization implies that
$\fo$ is equal to the minimal cardinality of a $\tau$-family that is not $o$-diagonalizable.
A detailed study of $\fo$ is initiated; the main results being that
consistently $\fo<\min\{\h,\s,\b\}$, and that in many standard models of set theory
(Cohen, Random, Hechler, Laver, Mathias, and Miller), $\cov(\M)=\fo$.

\begin{prob}[Problem 4.13 in \upannouncement{cot0}]\label{covMo}
Is $\cov(\M)=\fo$?
\end{prob}

It is worthwhile mentioning that this problem, which originated from the topological
studies of the minimal tower problem, is of similar flavor: It is well-known that
if $\p=\aleph_1$, then $\t=\aleph_1$ too. We have a similar assertion for
$\cov(\M)$ and $\fo$: If $\cov(\M)=\aleph_1<\b$, then $\cov(\M)=\fo$.

\begin{table}[!ht]
\begin{changemargin}{-3cm}{-3cm}
\begin{center}
{\tiny
\begin{tabular}{|r||cccccccccccccccccccccc|}
\hline
   & \smb{0} & \smb{1} & \smb{2} & \smb{3} & \smb{4} & \smb{5} & \smb{6} & \smb{7} &
   \smb{8} & \smb{9} & \smb{10} & \smb{11} & \smb{12} & \smb{13} & \smb{14} & \smb{15} &
   \smb{16} & \smb{17} & \smb{18} & \smb{19} & \smb{20} & \smb{21}\cr
\hline\hline

\mb{ 0} &
\yup&\yup&\yup&\yup&\nop&\fbn&\fbn&\fbn&\nop&\nop&\nop&
\nop&\yup&\yup&\fbn&\mbq&\nop&\nop&\yup&\yup&\yup&\yup\\
\mb{ 1} &
\mbq&\yup&\yup&\yup&\nop&\fbn&\fbn&\fbn&\nop&\nop&\nop&
\nop&\yup&\yup&\fbn&\mbq&\nop&\nop&\mbq&\yup&\yup&\yup\\
\mb{ 2} &
\nop&\nop&\yup&\yup&\nop&\nop&\fbn&\fbn&\nop&\nop&\nop&
\nop&\nop&\yup&\nop&\mbq&\nop&\nop&\nop&\nop&\yup&\yup\\
\mb{ 3} &
\nop&\nop&\nop&\yup&\nop&\nop&\nop&\fbn&\nop&\nop&\nop&
\nop&\nop&\nop&\nop&\nop&\nop&\nop&\nop&\nop&\nop&\yup\\
\mb{ 4} &
\yup&\yup&\yup&\yup&\yup&\yup&\yup&\yup&\nop&\nop&\mbq&
\mbq&\yup&\yup&\yup&\yup&\nop&\mbq&\yup&\yup&\yup&\yup\\
\mb{ 5} &
\mbq&\yup&\yup&\yup&\mbq&\yup&\yup&\yup&\nop&\nop&\mbq&
\mbq&\yup&\yup&\yup&\yup&\nop&\mbq&\mbq&\yup&\yup&\yup\\
\mb{ 6} &
\nop&\nop&\yup&\yup&\nop&\nop&\yup&\yup&\nop&\nop&\mbq&
\mbq&\nop&\yup&\nop&\yup&\nop&\mbq&\nop&\nop&\yup&\yup\\
\mb{ 7} &
\nop&\nop&\nop&\yup&\nop&\nop&\nop&\yup&\nop&\nop&\nop&
\mbq&\nop&\nop&\nop&\nop&\nop&\nop&\nop&\nop&\nop&\yup\\
\mb{ 8} &
\yup&\yup&\yup&\yup&\yup&\yup&\yup&\yup&\yup&\yup&\yup&
\yup&\yup&\yup&\yup&\yup&\yup&\yup&\yup&\yup&\yup&\yup\\
\mb{ 9} &
\mbq&\yup&\yup&\yup&\mbq&\yup&\yup&\yup&\mbq&\yup&\yup&
\yup&\yup&\yup&\yup&\yup&\yup&\yup&\mbq&\yup&\yup&\yup\\
\mb{10} &
\nop&\nop&\yup&\yup&\nop&\nop&\yup&\yup&\nop&\nop&\yup&
\yup&\nop&\yup&\nop&\yup&\nop&\yup&\nop&\nop&\yup&\yup\\
\mb{11} &
\nop&\nop&\nop&\yup&\nop&\nop&\nop&\yup&\nop&\nop&\nop&
\yup&\nop&\nop&\nop&\nop&\nop&\nop&\nop&\nop&\nop&\yup\\
\mb{12} &
\mbq&\mbq&\mbq&\mbq&\nop&\fbn&\fbn&\fbn&\nop&\nop&\nop&
\nop&\yup&\yup&\fbn&\mbq&\nop&\nop&\mbq&\yup&\yup&\yup\\
\mb{13} &
\nop&\nop&\nop&\nop&\nop&\nop&\nop&\nop&\nop&\nop&\nop&
\nop&\nop&\yup&\nop&\mbq&\nop&\nop&\nop&\nop&\yup&\yup\\
\mb{14} &
\mbq&\mbq&\mbq&\mbq&\fbn&\fbn&\fbn&\fbn&\nop&\nop&\fbn&
\fbn&\yup&\yup&\yup&\yup&\nop&\mbq&\mbq&\yup&\yup&\yup\\
\mb{15} &
\nop&\nop&\nop&\nop&\nop&\nop&\nop&\nop&\nop&\nop&\nop&
\nop&\nop&\yup&\nop&\yup&\nop&\mbq&\nop&\nop&\yup&\yup\\
\mb{16} &
\mbq&\mbq&\mbq&\mbq&\mbq&\mbq&\mbq&\mbq&\mbq&\mbq&\mbq&
\mbq&\yup&\yup&\yup&\yup&\yup&\yup&\mbq&\yup&\yup&\yup\\
\mb{17} &
\nop&\nop&\nop&\nop&\nop&\nop&\nop&\nop&\nop&\nop&\nop&
\nop&\nop&\yup&\nop&\yup&\nop&\yup&\nop&\nop&\yup&\yup\\
\mb{18} &
\nop&\nop&\nop&\nop&\nop&\nop&\nop&\nop&\nop&\nop&\nop&
\nop&\nop&\mbq&\nop&\mbq&\nop&\nop&\yup&\yup&\yup&\yup\\
\mb{19} &
\nop&\nop&\nop&\nop&\nop&\nop&\nop&\nop&\nop&\nop&\nop&
\nop&\nop&\mbq&\nop&\mbq&\nop&\nop&\nop&\yup&\yup&\yup\\
\mb{20} &
\nop&\nop&\nop&\nop&\nop&\nop&\nop&\nop&\nop&\nop&\nop&
\nop&\nop&\mbq&\nop&\mbq&\nop&\nop&\nop&\nop&\yup&\yup\\
\mb{21} &
\nop&\nop&\nop&\nop&\nop&\nop&\nop&\nop&\nop&\nop&\nop&
\nop&\nop&\nop&\nop&\nop&\nop&\nop&\nop&\nop&\nop&\yup\\

\hline
\end{tabular}
}
\end{center}
\end{changemargin}
\caption{Known implications and nonimplications}\label{imptab}
\end{table}

\subsection{A modified table of open problems}
It is possible that the diagram in Figure \ref{tauSch} is incomplete:
There are many unsettled possible implications in it.
After \cite{tautau, ShTb768}, there remained $76$ potential implications
which were not proved or ruled out.
The mentioned results of \upannouncement{cot0}
rule out $21$ of these implications (the reasoning is as follows:
If $P$ and $Q$ are properties with $\non(P)<\non(Q)$ consistent,
then $Q$ does not imply $P$), so that $55$ implications remain
unsettled.
The situation is summarized in Table \ref{imptab}, which updates the corresponding
table given in Issue 4 of the \SPMBul{}.

Each entry $(i,j)$ ($i$th row, $j$th column) contains a symbol.
\checkmark means that property $(i)$ in Figure \ref{tauSch} implies
property $(j)$ in Figure \ref{tauSch}.
$\times$ means that property $(i)$ does not (provably) imply property $(j)$,
and \textbf{?} means that the corresponding implication is still unsettled.
The new results are framed.

\begin{prob}
Settle any of the remaining $55$ implication in Table \ref{imptab}.
\end{prob}

\section{Problem of the month}

See Problem \ref{covMo} above.

\section{Problems from earlier issues}
In this section we list the past problems posed in the \SPMBul{},
in the section \probmonth{}.
For definitions, motivation and related results, consult the
corresponding issue.

For conciseness, we make the convention that
all spaces in question are
zero-dimentional, separable metrizble spaces.

\begin{issue}
Is $\binom{\Omega}{\Gamma}=\binom{\Omega}{\Tau}$?
\end{issue}

\begin{issue}
Is $\ufin(\Gamma,\Omega)=\sfin(\Gamma,\Omega)$?
And if not, does $\ufin(\Gamma,\Gamma)$ imply
$\sfin(\Gamma,\Omega)$?
\end{issue}

\begin{issue}
Does there exist (in ZFC) a set satisfying
$\ufin(\O,\O)$ but not $\ufin(\O,\Gamma)$?
\end{issue}
\begin{proof}[Solution]
\textbf{Yes} (Lubomyr Zdomsky).
\end{proof}

\begin{issue}
Does $\sone(\Omega,\Tau)$ imply $\ufin(\Gamma,\Gamma)$?
\end{issue}

\begin{issue}
Is $\p=\p^*$?
\end{issue}

\begin{issue}
Does there exist (in ZFC) an uncountable set satisfying $\sone(\BG,\B)$?
\end{issue}

\begin{issue}
Assume that $X$ has strong measure zero and $|X|<\b$.
Must all finite powers of $X$ have strong measure zero?
\end{issue}
\begin{proof}[Solution]
\textbf{Yes} (Scheepers; Bartoszy\'nski).
\end{proof}

\begin{issue}
Does $X \nin \NON(\M)$ and $Y\nin\mathsf{D}$ imply that
$X\cup Y\nin \COF(\M)$?
\end{issue}

\begin{issue}
Is $\split(\Lambda,\Lambda)$ preserved under taking finite unions?
\end{issue}
\begin{proof}[Solution]
Consistently yes (Zdomsky \upannouncement{addsplit}).
We conjecture that it is also consistently no. (This should hold under CH.)
\end{proof}

\general

\end{document}